\documentclass{lmcs} 
\pdfoutput=1

\usepackage{lastpage}
\lmcsdoi{15}{2}{7}
\lmcsheading{}{\pageref{LastPage}}{}{}%
{Aug.~04,~2017}{Apr.~30,~2019}{}

\pdfoutput=1


\usepackage{microtype}
\usepackage{hyperref}

\def\M{{\EM{\mathcal{M}}}}
\def\N{{\EM{\mathcal{N}}}}

\def\rmL{{\EM{\mathrm{L}}}}
\def\rmV{{\EM{\mathrm{V}}}}

\def\cS{{\EM{\mc{S}}}}

\def\aitch{{\EM{\downarrow}}}
\def\nh{{\EM{\uparrow}}}

\newcommand{\brac}[1]{[#1]}

\newcommand{\hatalpha}{{\widehat{\alpha}}}
\newcommand{\hatbeta}{{\widehat{\beta}}}
\newcommand{\hatgamma}{{\widehat{\gamma}}}

\newcommand{\incode}{\mathbf{in}}
\newcommand{\borelcode}{\mathbf{bor}}

\newcommand{\unioncode}{\mathbf{cup}}
\newcommand{\negcode}{\mathbf{neg}}
\newcommand{\ittmcode}{\mathbf{ittm}}
\newcommand{\id}{\mathrm{id}}

\newcommand{\pars}{\,\cdot\,}

\newcommand{\defn}[1]{{\bf{#1}}}
\newcommand{\defas}{{\EM{:=}}}

\newcommand{\BC}{\mathrm{BC}}

\newcommand{\ITTM}{\mathrm{ITTM}}
\newcommand{\LEFT}{\mathrm{LEFT}}
\newcommand{\RIGHT}{\mathrm{RIGHT}}
\newcommand{\STAY}{\mathrm{STAY}}

\newcommand{\Pfin}{\mathcal{P}_{<\omega}}

\def\Cantor{\EM{2^{\w}}}

\def\st{\EM{:}}

\def\nl{\newline}

\def\ORD{\mathrm{ORD}}

\def\dom{\text{dom}}

\def\realize{{\EM{{\mbf{R}}}}}

\DeclareMathOperator{\rank}{rank}

\def\w{\EM{\omega}}

\def\Naturals{{\EM{{\mbb{N}}}}}
\def\Nats{\Naturals}

\def\^{\EM{{}^{\And}}}

\def\And{\EM{\wedge}}

\def\<{\EM{\langle}}
\def\>{\EM{\rangle}}

\def\indent{\hspace*{2em}}

\def\EM#1{\ensuremath{#1}}
\def\mbb#1{\EM{\mathbb{#1}}}
\def\mbf#1{\EM{\mathop{\pmb{#1}}}}
\def\mc#1{\EM{\mathcal{#1}}}

\def\ul#1{\underline{#1}}

\begin{document}

\title[Feedback computability on Cantor space]{Feedback computability on Cantor space}

\author[Ackerman]{Nathanael L.\ Ackerman\rsuper{a}}
\author[Freer]{Cameron E.\ Freer\rsuper{b}}
\author[Lubarsky]{Robert S.\ Lubarsky\rsuper{c}}
\address{\lsuper{a}Department of Mathematics\\
Harvard University\\
Cambridge, MA 02138, USA}
\email{nate@math.harvard.edu}

\address{\lsuper{b}Department of Brain and Cognitive Sciences\\
Massachusetts Institute of Technology\\
Cambridge, MA 02139, USA}
\email{freer@mit.edu}

\address{\lsuper{c}Department of Mathematical Sciences\\
Florida Atlantic University\\
Boca Raton, FL 33431, USA}
\email{Robert.Lubarsky@alum.mit.edu}




\begin{abstract}
  \noindent 
	We introduce the notion of feedback computable functions from $\Cantor$ to  $\Cantor$, extending feedback Turing computation in analogy with the standard notion of computability for functions from $\Cantor$ to $\Cantor$.
	We then show that the feedback computable functions are precisely the effectively Borel functions. 	With this as motivation we  define the notion of a feedback computable function on a structure, independent of any coding of the structure as a real.  We show that this notion is absolute, and as an example characterize those functions that are computable from a Gandy ordinal with some finite subset distinguished.
\end{abstract}

\maketitle



\section{Feedback machines and Borel maps}

One of the most important observations of (effective) descriptive
set theory  is that every continuous map between Polish spaces is
computable with respect to some oracle, and every map which is
computable with respect to some oracle is continuous. (See \cite[Ex.\ 3D.21]{M}.)

This fact allows one to transport results from computability
theory to  the theory of continuous functions, and vice versa. But
it also is important because it provides a machine model for
continuity. Specifically, it provides a way of thinking about a continuous
map on a Polish space as a construction of the output from the
input, instead of simply via the abstract definition (requiring the
inverse image of open sets to be open).

We will show that this correspondence extends to feedback Turing computation \cite{AFL} and Borel functions. That is, we will show that feedback Turing computability provides a machine model for Borel functions on $\Cantor$, and conversely that every feedback computable function is itself Borel. These results should not be surprising, as it is already known (as reviewed below) that feedback computable reals are exactly the hyperarithmetic (i.e., $\Delta^1_1$) reals, and it is an old result of descriptive set theory that the Borel sets are exactly the $\Delta^1_1$-definable sets. So what we are doing here could be viewed as merely a type shift --- from $\Delta_1^1$ reals, which are $\Delta_1^1$ properties of natural numbers, to $\Delta^1_1$ properties on reals (for sets, or on pairs of reals if thinking about the graph of a function).

\subsection{Notation}

It will be useful to let $\tau\colon\Nats \times \Nats \to \Nats$ be a computable bijection. We will also need a computable bijection $\iota\colon\Nats \to 2^{<\w}$. For $\sigma \in 2^{<\w}$, define $\brac{\sigma} \defas \{x \in \Cantor \st \sigma \prec x\}$, i.e., the collection of elements of $\Cantor$ extending $\sigma$.

For a set $A$, define $\Pfin(A)$ to be the finite powerset of $A$, i.e., the set of all finite subsets of $A$.

By a \defn{countable ordinal} we will
mean a well-founded linear order $(A, \lhd)$ such that $A
\subseteq \w$. In particular, this representation of countable
ordinals will make it possible for oracle machines, as well as
feedback machines, to access them as oracles.

Given an admissible ordinal $\alpha$, let $\alpha^+$ be the \defn{next admissible after $\alpha$}, i.e., the least admissible ordinal greater than $\alpha$. An ordinal $\alpha$ is defined to be a \defn{Gandy ordinal} \cite{MR0432431} if it is admissible and for all $\gamma < \alpha^+$ there is an $\alpha$-computable well-ordering of order type $\gamma$.

For more details on admissible sets and ordinals, and on (effective) descriptive set theory, see \cite{HRT}, \cite{Kechris}, and \cite{M}.

\subsection{Feedback computability}

We now review the notion of feedback computability studied in \cite{AFL}. The intuitive idea is that we want to make sense of the notion of a machine which can ask halting queries of \emph{machines of the same type}.

The notation
$\{e\}^X_F(n)$ denotes the $e$th Turing machine with oracle $X$ and
halting function $F$ (which can also be interpreted as an oracle) on input $n$. When $\{e\}_F^X(n)$ queries $X$ it is said to be making an \defn{oracle query}, and when it queries $F$ it is said to be making a \defn{halting query}.

\begin{defi}
For any
$X\colon\w \rightarrow \{0, 1\}$ define the set $H_X \subseteq
\w \times \w$ to be the smallest collection for which there is a function $h_X\colon H_X \rightarrow
\{\uparrow,\downarrow\}$ satisfying the following:

\begin{enumerate}

\item[($\downarrow$)] If $\{e\}^X_{h_X}(n)$ makes no halting queries outside of $H_X$ and
converges after a finite number of steps then $(e,n) \in H_X$ and
$h_X(e,n) = \downarrow$, and conversely.

\item[($\uparrow$)] If $\{e\}^X_{h_X}(n)$ makes no halting queries outside of $H_X$ and
does not converge (i.e., runs forever) then $(e,n) \in H_X$ and
$h_X(e,n) = \uparrow$, and conversely.
\end{enumerate}
Furthermore, this $h_X$ is unique.
\end{defi}

\begin{defi}
A \defn{feedback Turing machine} (or \defn{feedback machine} for
short) is a machine of the form $\{e\}^X_{h_X}$ for some $e \in
\w$. The notation $\<e\>^X(n)$ is shorthand for
$\{e\}^X_{h_X}(n)$.

The set $H_X$ is the collection of \defn{non-freezing} computations
and the notation $\<e\>^X(n)\Downarrow$ means $(e,n) \in H_X$. If
$(e,n) \not \in H_X$ then $\<e\>^X(n)$ is
\defn{freezing}, written $\<e\>^X(n)\Uparrow$.

\end{defi}

One of the most important results about feedback computability is that feedback reducibility is equivalent to $\Delta^1_1$- or hyperarithmetic reducibility.   A set $X$ is hyperarithmetically reducible to $Y$ if $X$ is in the least admissible set containing $Y$.

\begin{thm}[{\cite[Theorem~16]{AFL}}]
\label{feedback=delta11}
For any $X, Y\colon\w \to \{0, 1\}$ the following are equivalent.
\begin{itemize}
\item $X$ is $\Delta^1_1(Y)$.

\item There is a feedback machine $e$ such that $\<e\>^Y(n) = X(n)$ for all $n \in \Nats$.
\end{itemize}
\end{thm}

Having recalled the notion of feedback computability, we now introduce the notion of a feedback computable function from $\Cantor$ to $\Cantor$.

\begin{defi}
	For $X\subseteq \Nats$,
a map $f\colon \Cantor \to \Cantor$ is \defn{feedback computable with respect
to $X$} if there is a feedback machine with
code $e$ such that for all $Y \in \Cantor$, the function $\<e\>^{X, Y}$ is total and for all $n \in \w$,  $\<e\>^{X, Y}(n)  = f(Y)(n)$. In this case $e$ is said to \defn{code} (with respect to $X$) a feedback computable function from $\Cantor$ to $\Cantor$.
\end{defi}

In other words, a function is feedback computable if there is a
feedback machine which, when given a description of a point in the
domain, outputs a description of the image of the point in the
range.

\subsection{Borel codes}

A \emph{Borel code} of a Borel subset of $\Cantor$ captures the way in which the Borel set was built up from basic open sets using the operations of countable union and complementation. There are many different types of Borel codes, all of which are, for practical purposes, equivalent. However, for our purposes it will be convenient to give a concrete coding system where each Borel code is a an element of $\Nats^\Nats$.

We begin with a couple of basic operations on functions from $\Nats$ to $\Nats$. Suppose $f\colon \Nats \to \Nats$. Let $f_*\colon \Nats \to \Nats$ be such that
$f_*(n) \defas f(n+1)$ for all $n \in \Nats$. Also, for $m \in \Nats$, let $f_m\colon \Nats\to \Nats$ be such that $f_m(n) \defas f(\tau(m, n) +1)$ for all $n \in \Nats$.

\begin{defi}
Let $\BC\subseteq \Nats^\Nats$ be
	the collection of \defn{Borel codes} (for $\Cantor$),
	defined by induction as follows.

\begin{itemize}
\item $\BC_0 \defas \{f \st f(0) \geq 2\}$.

\item If $\beta \in \ORD$ is greater than $0$ then $\BC_{\beta}$ is the smallest set containing
\begin{itemize}
\item $B_{\beta^*} \cup \{f \st f(0) =  0 \And f_* \in \BC_{\beta^*}\}$ for $\beta^* < \beta$, and

\item all functions $f\colon \Nats \to \Nats$ such that $f(0) = 1$ and for all $m \in \Nats$ there is a $\beta_n < \beta$ such that $f_m \in B_{\beta_n}$.
\end{itemize}
\end{itemize}

Finally, set $\BC \defas \bigcup_{\alpha < \w_1} \BC_\alpha$. If $f \in \BC$ then define the \defn{rank} of $f$ to be the least ordinal $\alpha$ such that $f \in \BC_\alpha$.
\end{defi}

A Borel code (for $\Cantor$) has an associated Borel set, called its \emph{realization}.

\begin{defi}
Suppose $\zeta\in \Nats^\Nats$ is a Borel code (for $\Cantor$). Define the \defn{realization} of $\zeta$, written $\realize(\zeta)$, to be the Borel subset of $\Cantor$
defined by induction on the rank of the $\zeta$ as follows.
\begin{itemize}
\item If $\zeta(0) \geq 2$,
		then $\realize(\zeta) \defas \brac{\sigma}$ where $\sigma = \iota(\zeta(0) - 2)$. Note that in this case,
	$\rank(\zeta) = 0$.

\item If $\zeta(0) = 0$,
		then let $\realize(\zeta) \defas \Cantor \setminus \realize(\zeta_*)$. Note that
in this case, if $\rank(\zeta) = \beta+1$ for some $\beta \in \ORD$,
then		$\rank(\zeta_*) = \beta$.

\item If $\zeta(0) = 1$ then $\realize(\zeta) \defas \bigcup_{m \in \Nats} \realize(\zeta_m)$. Note that in this case, $\rank(\zeta_m) < \rank(\zeta)$ for all $m \in \Nats$.

\end{itemize}
\end{defi}

By the \emph{union} or \emph{intersection} of a collection of Borel codes, we will mean the code for union or intersection of their realizations.

One of the first steps in showing that the Borel functions and the feedback computable functions coincide is to show that feedback computability interacts well with Borel codes.

For any ordinal $\alpha\in\w_1$, define an \defn{encoding}
of $\alpha$ to be a linear ordering on a subset of $\w$ of order type $\alpha$.

\begin{lem}
	\label{borel-code-for-ordinal}
There is a feedback machine $\borelcode$ such that for any $\alpha \in \w_1$,
	any encoding $\hatalpha$ of $\alpha$,
	any $X\in \Cantor$, and any $n \in \Nats$ such that $\<n\>^X$ is total, we have
\begin{itemize}
\item $\<\borelcode\>^{X, \hatalpha}(n) = 1$ if $\<n\>^X \in \BC_{\alpha}$, and

\item $\<\borelcode\>^{X, \hatalpha}(n) = 0$ if $\<n\>^X \not\in \BC_{\alpha}$.

\end{itemize}
\end{lem}
\begin{proof}
Let $\<\borelcode\>^{X, \hatalpha}(n)$ do the following. \nl\nl
\ul{Step 1:} \nl
If $\<n\>^X(0) \geq 2$ then return $1$. \nl\nl
\ul{Step 2:} \nl
If $\<n\>^X(0) = 0$ then search for a $\beta < \alpha$ such that $\<\borelcode\>^{X, \hatbeta}(n_*) = 1$ where
$\hatbeta$ is the encoding of $\beta$ induced by $\hatalpha$, and where
	$\<n_*\>^X(m) \defas \<n\>^X(m+1)$ for all $m \in \Nats$. If there exists such a $\beta$ then return $1$, and otherwise return $0$. \nl\nl
\ul{Step 3:} \nl
	If $\<n\>^X(0) =1$ then for each $m \in \Nats$ search for a $\beta_m < \alpha$ such that $\<\borelcode\>^{X, \widehat{\beta_m}}(n_m) = 1$ where
	$\widehat{\beta_m}$ is the encoding of $\beta_m$ induced by $\hatalpha$, and where $\<n_m\>^X(k) \defas \<n\>^X(\tau(m, k) + 1)$ for all $m \in \Nats$.
If there exists such a $\beta_m$ for each $m \in \Nats$ then return $1$, and otherwise return $0$. \nl\nl
It is an easy induction to show that $\borelcode$ is the desired code.
\end{proof}

In other words, there is a feedback machine which, uniformly in $\alpha$ and an oracle $X$, can check whether or not a total function feedback computable in $X$ is a Borel code of rank at most $\alpha$.

\begin{lem}
\label{Element of Borel code is feedback computable}
There is a feedback machine $\incode$ such that for any Borel code $C$, any $X\in \Cantor$ and any $n \in \Nats$ such that $\<n\>^X$ is total, we have
\begin{itemize}
\item $\<\incode\>^{C, X}(n) = 1$ if $\<n\>^X \in \realize(C)$, and

\item $\<\incode\>^{C, X}(n) = 0$ if $\<n\>^X \not\in \realize(C)$.

\end{itemize}
\end{lem}
\begin{proof}
Call $C$ the first oracle and $X$ the second oracle. We will need to unravel the code $C$. However, there is no mechanism for changing an oracle. This ends up not being a problem, because the changes we would like to make are simple --- which we formalize using the Recursion Theorem.

Toward this end, suppose we have computer code $e$, which we will think of as instructions for a feedback machine $\langle e \rangle$. We will define new code $e^*$, as follows. The behavior of the computation $\<e^*\>^{C,X}$ depends on the value of $C(0)$. \nl\nl
\ul{Case 1:} $C(0) \geq 2$.\nl
Let $\sigma \defas \iota(C(0) - 2)$. If $\<n\>^X(m) = \sigma(m)$ for all $m \in \dom(\sigma)$, then return $1$. Otherwise return $0$.
\nl\nl
\ul{Case 2:} $C(0) = 0$.\nl
Let $\<e\>^{C_*,X}$ be the machine that runs just like $\<e\>$ with oracle $C,X$,  except that whenever $e$ makes a query of $k$ to the first oracle, a query  of $k+1$ is made instead. Return $1 - \<e\>^{C_*,X}$.
\nl\nl
\ul{Case 3:} $C(0) = 1$.\nl
Let $\<e\>^{C_m,X}$ be the machine that runs just like $\<e\>$ with oracle $C,X$,  except that whenever $e$ makes a query of $k$ to the first oracle, a query  of $\tau(m,k)+1$ is made instead. Let $\<e^*\>^{C,X}$ search for an $m$ such that $\<e\>^{C_m,X} = 1$, and if it finds such an $m$ it returns 1, else 0.
\nl\nl
The function from $e$ to $e^*$ is computable. Let $\incode$ be a fixed point. It is an easy induction on the rank the Borel code $C$ to show that $\incode$ has the desired properties.
\end{proof}

In other words, there is a feedback machine which can determine whether or not a feedback computable function (relative to an oracle) is in the realization of a Borel code, uniformly in the Borel code and the oracle. Now we show, unsurprisingly, that if given a feedback computable collection of Borel codes, they can be combined to form a new Borel code.

\begin{lem}
There are feedback machines $\negcode, \unioncode$ such that
\begin{itemize}
\item if $\<n\>^X \in \BC$, then $\<\negcode\>^X(n, \pars) \in \BC$ with $\realize\bigl(\<\negcode\>^X(n, \pars)\bigr) = \Cantor \setminus \realize\bigl(\<n\>^X\bigr)$, and

\item if
	$\<n\>^X(m, \pars) \in \BC$
	for all $m \in \Nats$,
		then $\<\unioncode\>^X(n, \pars) \in \BC$ with $\realize\bigl(\<\unioncode\>^X(n, \pars)\bigr) = \bigcup_{m \in \w} \realize\bigl(\<n\>^X(m, \pars)\bigr)$.
\end{itemize}
\end{lem}

\begin{proof}
Let $\<\negcode\>^X(n, 0) = 0$ and $\<\negcode\>^X(n, m+1) = \<n\>^X(m)$ for all $m\in \Nats$. Let $\<\unioncode\>^X(n, 0) = 1$ and $\<\unioncode\>^X(n, m+1) = \<n\>^X(\tau^{-1}(m))$.
\end{proof}

We now define the notion of $\Delta^1_1(X)$ function. We do this in terms of Borel codes, since the Borel subsets of $\Cantor$ are exactly all the $\Delta^1_1$-sets. For more on this notion of $\Delta_1^1$-function, see \cite[Section~3D and Theorem~3E.5]{M}.

\begin{defi}
A map $f\colon \Cantor \to \Cantor$ is \defn{$\Delta^1_1(X)$} for $X \subseteq \Nats$ if there is a $\Delta^1_1(X)$ sequence of functions $(\gamma_\sigma)_{\sigma\in 2^{<\w}}$ such that $\realize(\gamma_\sigma) = f^{-1}(\brac{\sigma})$.
\end{defi}

The following result is standard (see, e.g., \cite[Chapter~2]{M}).

\begin{lem}
A map $f\colon \Cantor \to \Cantor$ is Borel if and only if $f$ is $\Delta^1_1(X)$ for some $X \subseteq \Nats$.
\end{lem}

\subsection{A characterization of feedback computable functions}

We now show that we can isolate Borel codes for those oracles that cause a feedback computation to halt with a tree of subcomputations of height at most $\alpha$.

\begin{prop}
\label{Borel codes of oracles satisfying a condition}
There is a computable
	collection of codes for feedback machines, $e_{\aitch}$,
$e_{\aitch j}$, $e_{\nh}$, $e_{\nh j}$ for $j \in \w$, such that
\begin{itemize}
\item for all $\alpha\in\w_1$ and encodings $\hatalpha$,

\item for all $X \subseteq \w$,

\item for all $\sigma \in 2^{<\w}$,

\item for all $f \in \w$, and

\item for all $x \in \{\aitch, \nh\}\cup \{\aitch, \nh\}\times \w$,
\end{itemize}
$\<e_x\>^{X, \hatalpha}(f,n, \sigma, \pars)$ is a Borel code, which we will refer to as $\zeta_x$.
	
	Further, whenever $Y \in \Cantor$ with
$\sigma \prec Y$ and $j\in\w$  the following properties hold.
\begin{itemize}
\item $Y \in \realize(\zeta_\aitch)$ if and only if $\<f\>^{X, Y}(n)$ halts with a tree of subcomputations of height $\leq \alpha$.

\item $Y \in \realize(\zeta_{\aitch j})$ if and only if $\<f\>^{X, Y}(n)$ halts with a tree of subcomputations of height $\leq \alpha$ and outputs $j$.

\item $Y \in \realize(\zeta_{\nh})$ if and only if $\<f\>^{X, Y}(n)$ does not halt and the tree of subcomputations is of height $\leq \alpha$.

\item $Y \in \realize(\zeta_{\nh j})$ if and only if $\<f\>^{X, Y}(n)$ does not halt after $j$-many steps, and up to the $j$\emph{th} step the tree of subcomputations has height $\leq \alpha$.
\end{itemize}

	Finally, we always have $Y \not \in \realize(\zeta_x)$ if $\sigma \not \prec Y$.
\end{prop}
\begin{proof}
We begin with some notation. Enumerate all triples $(\eta, \nu, k)$ such
that $\eta$ and $\nu$ are
elements of $2^{<\w}$ and $\{f\}^{X, \eta}_{\nu}(n)$ does not make any invalid oracle
calls in the first $k$-many steps, i.e., does not make any oracle queries outside of $\eta$ or $\nu$. Call this collection $B$.

For each
$\mbf{\eta} = (\eta, \nu, k) \in B$, let $\<(r_i^{\mbf{\eta}},
m_i^{\mbf{\eta}})\>_{i \leq \ell^{\mbf{\eta}}}$ be the sequence such that $\<\tau(r_i^{\mbf{\eta}},
m_i^{\mbf{\eta}})\>_{i \leq \ell^{\mbf{\eta}}}$ are the queries made
by $\{f\}^{X, \eta}_{\nu}(n)$ to $\nu$. Note that the length of the sequence is $\ell^{\mbf{\eta}}$.  Recall these are called halting queries. In particular if $\nu$ is the correct
response to every halting query made by $\<f\>^{X, Y}(n)$, i.e., $\nu$ agrees with $h_{X, Y}$, then
$\{f\}^{X, Y}_{\nu}(n) \cong \<f\>^{X, Y}(n)$. With this notation we will use the convention that $\nu(\tau(a,b)) = 0$ means $\nu$ ``believes'' $\<a\>^{X, Y}(b)$ halts and $\nu(\tau(a,b)) = 1$ means $\nu$ ``believes'' $\<a\>^{X, Y}(b)$ does not halt.

Our goal will be to define Borel codes $C_{\mbf{\eta}}^\hatalpha$ for all $\mbf{\eta} \in B$ in such a way that $C_{\mbf{\eta}}^\hatalpha$ is the Borel code of all $Y$ extending $\eta$ such that the behavior of halting calls of $\<f\>^{X, Y}(n)$ on the first $k$-many steps agrees with that of $\nu$, and for which the tree of computations has height at most $\alpha$.

Let $B_\aitch$ be the collection of triples $(\eta, \nu, k) \in B$ such that $\{f\}^{X, \eta}_{\nu}(n)$ halts in at most $k$ steps. For $j
\in \w$, let $B_{\aitch j}$ be the collection of triples in $B_\aitch$ with output $j$. Let $B_{\nh j}$ be the collection of
triples in $B\setminus B_\aitch$ whose third coordinate is $j$.

\ \\
\noindent For each $\mbf{\eta} = (\eta, \nu, k) \in B$ we define a Borel
code $C_{\mbf{\eta}}^\hatalpha$ as follows.
\begin{enumerate}
\item[(i)] If $\<f\>^{X, \eta}(n)$ makes no halting queries, then $C_{\mbf{\eta}}^\hatalpha$ is
a code for $\brac{\eta}$.

\item[(ii)] If $\{f\}^{X, \eta}_{\nu}$ makes any halting queries and $\alpha = 0$,
then $C_{\mbf{\eta}}^\hatalpha$ is a Borel code for $\emptyset$.

\item[(iii)] Suppose $\nu(\tau(r_i^{\mbf{\eta}}, m_i^{\mbf{\eta}})) = 0$, i.e., $\{f\}^{X, \eta}_{\nu}$
``thinks'' $\<r_i^{\mbf{\eta}}\>^{X, Y}(m_i^{\mbf{\eta}})$ halts. Then
let $D_{\mbf{\eta}, i}^\hatalpha$ be the Borel code which is the union of the Borel codes of
$\<e_\aitch\>^{X, \hatbeta}(r_i^{\mbf{\eta}}, m_i^{\mbf{\eta}}, \eta, \pars)$ for $\beta < \alpha$, where $\hatbeta$ is the encoding of $\beta$ induced by $\hatalpha$. In other words $D_{\mbf{\eta}, i}^\hatalpha$ is a Borel code for the collection of those $Y$ such that $\<r_i^{\mbf{\eta}}\>^{X, Y}(m_i^{\mbf{\eta}}) =0$ and has a tree of subcomputations of rank less than $\alpha$. Or said another way, $D_{\mbf{\eta}, i}^\hatalpha$ is a Borel code for the collection of those $Y$ such that $\<f\>^{X, Y}(n)$ agrees with $\nu$ on the $i$th halting query and has a tree of subcomputations of rank less than $\alpha$.

Suppose $\nu(\tau(r_i^{\mbf{\eta}}, m_i^{\mbf{\eta}})) = 1$, i.e., $\{f\}^{X, \eta}_{\nu}$ ``thinks''
$\<r_i^{\mbf{\eta}}\>^{X, Y}(m_i^{\mbf{\eta}})$ does not halt. Then
for $j \in \w$ and $\beta < \alpha$, let $E^\hatalpha_{\mbf{\eta}, i, j, \hatbeta}$ be the Borel code of the union of $\<e_{\nh j}\>^{X, \hatbeta}(r_i^{\mbf{\eta}}, m_i^{\mbf{\eta}}, \eta, \pars)$. In other words $E_{\mbf{\eta}, i, j,\beta}^\hatalpha$ is a Borel code for those $Y$ such that $\<r_i^{\mbf{\eta}}\>^{X, Y}(m_i^{\mbf{\eta}})$ hasn't halted by step $j$ and which has a tree of subcomputations of height at most $\beta$. Now for $\gamma < \alpha$ let $D_{\mbf{\eta}, i}^{\hatalpha, \hatgamma}$ be the Borel code of the intersection of $E_{\mbf{\eta}, i, j, \beta}^\alpha$ over $j \in \w$ and $\beta < \gamma$, where $\hatgamma$ is the encoding of $\gamma$ induced by $\hatalpha$. So $D_{\mbf{\eta}, i}^{\hatalpha, \hatgamma}$ is a Borel code for those $Y$ such that $\<r_i^{\mbf{\eta}}\>^{X, Y}(m_i^{\mbf{\eta}})$ doesn't halt and has a tree of subcomputations of height at most $\gamma$. Now let $D_{\mbf{\eta}, i}^\hatalpha$ be the union of $D_{\mbf{\eta}, i}^{\hatalpha, \hatgamma}$ over $\gamma < \alpha$. So $D_{\mbf{\eta}, i}^{\hatalpha}$ is a Borel code for those $Y$ such that $\<r_i^{\mbf{\eta}}\>^{X, Y}(m_i^{\mbf{\eta}})$ doesn't halt and has a tree of subcomputations of height at less than $\alpha$.

Let $C_{\mbf{\eta}}^\hatalpha$ be the intersection of $D_{\mbf{\eta}, i}^\hatalpha$ for $i \leq \ell^{\mbf{\eta}}$. Note that as $\ell^{\mbf{\eta}}$ is finite the rank of the tree of subcomputations of anything in $C_{\mbf{\eta}}^\hatalpha$ bounded by the supremum of the tree of subcomputations of $\<r_i^{\mbf{\eta}}\>^{X, Y}(m_i^{\mbf{\eta}})$ plus $1$. Therefore the rank of the tree of subcomputations of $\<f\>^{X, Y}(n)$ is at most $\alpha$ for all $Y \in C_{\mbf{\eta}}^\hatalpha$.
\end{enumerate}

Finally for $j\in \omega$ we define the codes as follows.
\begin{enumerate}
\item[(I)] $\<e_\aitch\>^{X, \hatalpha}(f, n, \eta, \pars)$ is the union of $C_{\mbf{\eta}}^\hatalpha$ such
that $\mbf{\eta} \in B_\aitch$.

\item[(II)] $\<e_{\aitch j}\>^{X, \hatalpha}(f, n, \eta, \pars)$ is the union of $C_{\mbf{\eta}}^\hatalpha$
such that $\mbf{\eta} \in B_{\aitch j}$.

\item[(III)] $\<e_{\nh j}\>^{X, \hatalpha}(f, n, \eta, \pars)$ is the union of $C_{\mbf{\eta}}^\hatalpha$ such
that $\mbf{\eta} \in B_{\nh j}$.

\item[(IV)] $\<e_{\nh}\>^{X, \hatalpha}(f, n, \eta, \pars)$ is the intersection of the codes
given by $\<e_{\nh j}\>^{X, \hatalpha}(f, n, \eta, \pars)$ for $j \in \w$.
\end{enumerate}

\begin{clm}
\label{Borel codes of oracles satisfying a condition:Claim 1}
$C_{\mbf{\eta}}^\hatalpha$ is a Borel code of the set of
all $Y$ extending $\eta$ such that the behavior of $\<f\>^{X, Y}(n)$ on the
first $k$-many steps agrees with that of $\nu$, and for which the
tree of computations has height at most $\hatalpha$.
\end{clm}
\begin{proof}
Our proof that these codes satisfy our theorem proceeds by induction on $\alpha$.
Notice by conditions (i) and (ii) that if $\alpha = 0$  then
$C_{\mbf{\eta}}^\hatalpha$ satisfies the claim.

But then by induction on $\alpha$ and condition (iii),
$D_{\mbf{\eta}, i}^\hatalpha$  is a Borel code for those $Y$ extending
$\eta$ for which the $i$th halting call agree with $\nu$ and the
tree of subcomputations has height $< \alpha$. This then implies that $C_{\mbf{\eta}}^\hatalpha$ satisfies the claim.
\end{proof}

It is then straightforward to check from Claim~\ref{Borel codes of oracles satisfying a condition:Claim 1} that these definitions satisfy the proposition.
\end{proof}

Before moving on to the main application of Proposition~\ref{Borel
codes of oracles satisfying a condition}, it is worth taking a
moment to highlight the importance of the ordinal $\alpha$. Specifically, if we did not have a uniform bound on the height of the tree of subcomputations we were considering, we might accidentally make a halting query which would cause our computation to freeze --- causing the entire construction to break. However, we show in Lemma~\ref{bound on feedback tree for feedback computable functions} that this is not an issue, as there will always be a single bound on all trees of subcomputations, which is itself feedback computable from $X$.

\begin{lem}
\label{bound on feedback tree for feedback computable functions}
Suppose $e$ is a code (with respect
to $X$) for a feedback computable function from $\Cantor$ to $\Cantor$.  Then there is a countable ordinal $\alpha$ and an encoding $\hatalpha$ that is feedback
computable in $X$ such that for every $Y\in \Cantor$ and $n \in \Nats$ the tree
of subcomputations of $\<e\>^{X, Y}(n)$ has height bounded by $\alpha$.
\end{lem}

\begin{proof}
Let $P(T)$ be the statement ``$(\exists Y\in \Cantor)(\exists n \in \Nats)$ such that $T$ is locally the tree of subcomputations for $\<e\>^{X, Y}(n)$''. Then $P(\pars)$ is a $\Sigma^1_1(X)$ predicate. However, for all $Y \in \Cantor$ and $n \in \Nats$,  the feedback computation $\<e\>^{X, Y}(n)$ does not freeze and hence its tree of subcomputations is well-founded. In particular this implies that $P(T)$ holds only if $T$ is a well-founded relation. 

Hence, by a relativized version of \cite[II.5.9]{HRT},
there is some ordinal $\alpha$ hyperarithmetic in $X$ such that $\alpha$ bounds the height of all $T$ satisfying $P$. Further, by \cite[Theorem~16]{AFL}, some encoding $\hatalpha$ is feedback computable in $X$ because $\alpha$ is hyperarithmetic in $X$. In particular this implies that $\alpha$ bounds the height of the tree of subcomputations of $\<e\>^{X, Y}(n)$ for all $Y \in \Cantor$ and $n \in \Nats$.
\end{proof}

\begin{prop}
\label{Feedback_implies_Borel_for_maps} Suppose $f\colon2^\w \to
2^\w$ is a feedback computable map (with respect to $X$).  Then
	$f$ is $\Delta^1_1(X)$.
\end{prop}
\begin{proof}
Let $e$ be a code (with respect to $X$) for the map $f$, and let $\alpha$ and $\hatalpha$ be as in
Lemma~\ref{bound on feedback tree for feedback computable functions}.
	By Proposition~\ref{Borel codes of oracles satisfying a condition},
	there is a uniformly computable (in $X$ and $\hatalpha$) collection
of Borel codes $\zeta_{\aitch j}$ such that $\realize(\zeta_{\aitch j}) =
\{Y\st \<e\>^{X, Y}(n) = j$ and $\<e\>^{X, Y}(n)$ has a tree of
subcomputations of height $<\alpha\}$. But then
$\realize(\zeta_{\aitch j}) = \{Y\st \<e\>^{X, Y}(n) = j\}$, as the
trees of subcomputations for feedback machines of the form
$\<e\>^{X, Y}(n)$ have height $< \alpha$. Hence there is a
collection of Borel codes $(\gamma_\sigma)_{\sigma \in 2^{<\w}}$,
uniformly feedback computable in $X$, such that $\realize(\gamma_\sigma) =
	f^{-1}(\brac{\sigma})$ for each $\sigma \in 2^{<\w}$. But then by \cite[Theorem~16]{AFL}, the functions $\gamma_\sigma\colon \Nats \to \Nats$ are $\Delta^1_1(X)$ uniformly in $\sigma$. Therefore $f$ is $\Delta^1_1(X)$.
\end{proof}

\begin{prop}
\label{Borel_implies_Feedback_for_maps} Suppose $f\colon2^\w \to
2^\w$ is $\Delta^1_1(X)$.
Then $f$ is	feedback computable (with respect to $X$).
\end{prop}
\begin{proof}
	There is a
sequence $\<\gamma_\sigma\>_{\sigma \in 2^{<\w}}$ which is in $\Delta^1_1(X)$ such that for $\sigma \in 2^{<\w}$ the real
$\gamma_\sigma$ is a Borel code for $f^{-1}(\brac{\sigma})$. By Theorem~\ref{feedback=delta11}, we therefore have that $\<\gamma_\sigma\>_{\sigma \in 2^{<\w}}$ is feedback computable from $X$.

	Then by Lemma~\ref{Element of Borel code is feedback computable} there is a feedback machine $e$ such that $\<e\>^{X, Y}(n) = 1$ if there is a $\sigma\in 2^{n+1}$ such that $\sigma(n) = 1$ and $Y \in \realize(\gamma_{\sigma})$, and $\<e\>^{X, Y}(n) = 0$ otherwise. In other words, $\<e\>^{X, Y}(n)$ is the value of $f(Y)(n)$. Hence $e$ is a code (relative to $X$) for $f$.
\end{proof}

At this point, we have accomplished our main purpose in this section, to provide a machine model for Borel functions from $\Cantor$ to $\Cantor$, by showing that the Borel functions are exactly the feedback computable functions. To round this out, we include some other characterizations of feedback computable functions.

If $f$ is feedback computable (mention of the parameter $X \subseteq \omega$ will be suppressed), then $f(Y)$ is in any admissible set containing $Y$. So $f$ is uniformly $\Sigma_1$ definable over all admissible sets, as $f(Y) = Z$ iff within any admissible set containing $Y$ there is a tree witnessing the computation of $f(Y)$ \cite[Proposition~4]{AFL}. In fact, $f$ can trivially be extended to a function on the entire universe $V$, by letting $f(Y)$ be the empty set whenever $Y$ is not a real. Conversely, suppose $f \colon  V \rightarrow V$ is uniformly $\Sigma_1$ definable over all admissible sets, and $f$ takes reals to reals. Then, as a function on reals, $f$ is $\Sigma^1_1$; namely, $f(Y) = Z$ iff there is a real coding an $\omega$-standard admissible set modeling $f(Y) = Z$. (In a little more detail, it costs nothing to say the model is $\omega$-standard, as you can insist that the members of the model's version of $\omega$ are given by the evens in their natural order. Furthermore, the ordinal standard part of an admissible set is itself admissible, so an $\omega$-standard admissible set containing $Y$, even if non-standard, will also contain $f(Y)$.) It is folklore that every $\Sigma^1_1$ function is $\Delta^1_1$, as $f(Y) \not = Z$ iff there is a $W$ such that $f(Y) = W$ and $W \not = Z$.

For another characterization of the functions in question, van de Wiele \cite{vdw} showed that a function is uniformly $\Sigma_1$ over all admissible sets iff it is $E$-recursive. This was later extended by Slaman \cite {Sl2} (see also \cite {Lub} for a different
proof\footnote{There is a minor mistake in the latter which can easily be corrected. Slaman's proof uses selection from the parameter $p$. It is mistakenly claimed in \cite{Lub} that selection from $p$ is not necessary. In fact, the construction in \cite{Lub} is perfectly good; it's just that the use of selection from $p$ in the construction was overlooked.}) to include hereditarily countable parameters. Slaman's result is that for a hereditarily countable parameter $p$, a function $f$ is uniformly $\Sigma_1(p)$ over all admissible sets iff $f$ is $ES_p$ recursive, where $ES_p$ recursion is $E$-recursion  augmented by selection from $p$, a schema first identified in \cite {Hoo} and further studied in \cite {Sl1}. In our case, the parameter is a real $X$; by Gandy Selection, selection from a real follows from the regular $E$-recursion schema \cite {Sl1}, so that $f$ is uniformly $\Sigma_1(X)$ over all admissible sets iff $f$ is $E$-recursive in $X$.

Summarizing the above, we have the following theorem.

\begin{thm}
\label{Feedback_is_Borel_(for_maps)}
 For any function $f\colon 2^\w \to
2^\w$ and any $X \subseteq \w$ the  following are equivalent.

\begin{enumerate}
\item[(a)] $f$ is $\Delta^1_1(X)$.

\item[(b)] $f$ is feedback computable with respect to $X$.

\item [(c)] $f$ can be extended to a function on $V$ which is uniformly  $\Sigma_1(X)$ definable over all admissible sets.

\item [(d)] $f$ can be extended to a function $E$-recursive in $X$.
\end{enumerate}
\end{thm}

\section{Feedback computability relative to a structure}
In this section, we extend the notion of an oracle from a set of natural numbers to a structure (up to isomorphism). If the structure is countable, it can be coded as a set of natural numbers, however this cannot be done if the structure is uncountable. As such, we want our definition of computation from a structure to ultimately be independent of any coding of our structure. This can be seen as a feedback analogue of Medvedev reducibility on isomorphism classes of structures. (For a survey of Muchnik and Medvedev degrees, see
\cite{MR2931672}.) We will make use of the fact (Theorem~\ref{Feedback_is_Borel_(for_maps)}) that
Borel functions can be thought of as those that are feedback computable from an oracle, and the fact that the isomorphism classes of countable structures are Borel (see \cite[Theorem~16.6]{Kechris}).
In fact, using Theorem~\ref{Feedback_is_Borel_(for_maps)}, one can rephrase the results in this section in terms of notions from effective descriptive set theory, instead of feedback computability, if so desired.

\begin{defi}
	A countable language $L_0$ is \defn{feedback computable} if the sets of relation, function, and constant symbols in $L_0$, and their arities, are uniformly feedback computable.
\end{defi}
	
	A particular feedback computable enumeration of such data gives rise to a natural encoding of each countable $L_0$-structure with underlying set $\Nats$, which we will sometimes call its \defn{$L_0$-encoding}, and often use implicitly. (For more details on such encodings, see, e.g., \cite{AntonioBook}.)

A structure is hereditarily countable when its underlying set and the sets of its relations, functions, and constants are hereditarily countable.
Note that a structure is isomorphic to a hereditarily countable one if and only if it is a countable structure in a countable language.
However, there are hereditarily countable structures that are uncountable in some admissible sets that contain them as an element.
Hence we may think of a structure being hereditarily countable in an admissible set as a measure of its computability with respect to the admissible set.

We now give two definitions of different kinds of structures that can be computed from a structure independent of any coding.

\subsection{Feedback computing expansions}

We now introduce the notion of feedback computing an expansion of a structure.

\begin{defi}
	For $j\in\{0,1\}$, let
$L_j$ be a feedback computable language,
and let $\M_j^*$ be an $L_j$-structure with underlying set $\Nats$ and natural $L_j$ encoding $i_j$.
	We say that $e\in\Nats$ \defn{feedback computes $\M_1^*$ from $\M_0^*$ using the natural encodings} if
	$\<e\>^{i_0(\M_0^*)}$ is a total function with $\<e\>^{i_0(\M_0^*)} = i_1(\M_1^*)$.
\end{defi}

We will later use the following lemma.

\begin{lem}
\label{remark-on-fb-natural-encoding}
	For $j\in\{0,1\}$, let
$L_j$ be a feedback computable language,
and let $\M_j^*$ be an $L_j$-structure with underlying set $\Nats$. The statement ``$e$ feedback computes $\M_1^*$ from $\M_0^*$ using the natural encodings'' is $\Sigma_1$ over any admissible set containing both $\M_j^*$'s, as is the statement ``there exists an $\M_1^*$ such that $e$ feedback computes $\M_1^*$ from $\M_0^*$ using the natural encodings''. Further, for every $\M_0^*$ and $e$, there is at most one $\M_1^*$ for which $e$ feedback computes $\M_1^*$ from $\M_0^*$ using the natural encodings.
\end{lem}

\begin{proof}
Each convergent feedback computation $\<e\>^{i_0(\M_0^*)}(n)$ has a witness to its convergence and value in any such admissible set. So if $\<e\>^{i_0(\M_0^*)}$ is total, by admissibility a set of witnesses to the convergence of each $\<e\>^{i_0(\M_0^*)}(n)$ can be formed, witnessing the totality of $\<e\>^{i_0(\M_0^*)}$. It is then arithmetic in the output of $\<e\>^{i_0(\M_0^*)}$ whether that function is the natural encoding of a structure, and what that structure would be, in particular whether it equals $\M_1^*$. So the $\Sigma_1$ sentence claimed to existence is the existence of a witness or a construction that $\<e\>^{i_0(\M_0^*)} = i_1(\M_1^*)$.
\end {proof}

\begin{defi}
	\label{fc-expansions-of-structures}
Let $L_0 \subseteq L_1$ be feedback computable languages, $\M_0$ be a countable $L_0$-structure, and $\M_1$ a countable $L_1$-structure which is an expansion of $\M_0$, i.e., $\M_1|_{L_0} = \M_0$. Then
	\defn{$e$ feedback computes the expansion $\M_1$ of oracle $\M_0$}, written
	$\<e\>^{\M_0} \asymp^+ \M_1$,
	if for every $L_0$-structure $\M_0^*$ with underlying set $\Nats$
	and isomorphism $f\colon \M_0^* \to \M_0$,
	for all $L_1$-structures $\M_1^*$ such that $e$ feedback computes $\M_1^*$ from $\M_0^*$ using the natural encodings,
	we have
	$\M_1^* |_{L_0} = \M_0^*$ and
	$f$ is an isomorphism from $\M_1^* \to \M_1$.
\end{defi}

As an example, consider the case where $L_0 = (\id, \times)$ and $L_1 = (\id, \times, (\cdot)^{-1})$ with $\id$ a constant, $\times$ a binary function, and $(\cdot)^{-1}$ a unary function. Suppose $\M_0$ is a group in the language $L_0$ and $\M_1$ is the same group but in the language $L_1$ (i.e., with the inverse function). Then, as we can feedback compute the inverse function when we are passed a group, there is a feedback machine which computes the expansion $\M_1$ of $\M_0$.

The following generalization of the
relativized L\'evy--Shoenfield Absoluteness Theorem \cite[Theorem~$36^\prime$]{MR1492987}
will be used in Propositions~\ref{Absolutness_of_asymp+} and \ref{Calculating-one-structure-from-another-absolute}.
A version of it is stated as Exercise~15.14 of \cite{MR1492987} (in only the first two editions);
we include its proof for completeness.

\begin{thm}[{\cite[Ex.~15.14]{MR1492987}}]
	\label{SL-lemma}
Let $A$ be a transitive admissible set containing the countable ordinals and in which $p$ is hereditarily countable. If $\varphi$ is a $\Sigma_1(p)$ sentence true in $V$, then $\varphi$ is true in $A$.
\end{thm}

\begin{proof}
Let $\varphi(p)$ be $\exists x \; \psi(x)$, where $\psi$ is $\Delta_0$ and in prenex form. Replace the existential quantifiers in $\psi$ by terms for Skolem functions. By L\"owenheim--Skolem, and using the hereditary countability of $p$, there is a witness $x$ to $\psi(x)$ in a countable transitive model. In particular, the height of this model is a countable ordinal, say $\gamma$. Work relative to some coding $Z$ of $p$, which by hypothesis exists in $A$. Build the tree of finite structures (allowing the Skolem functions to be partial), which we think of as finite substructures of models of $\psi(x)$, with a ranking function (of the sets of the model) into $\gamma$. There are two additional constraints on this tree. A node on level $n$ must contain (at least) the first $n$-many elements (according to $Z$) of the transitive closure of $p$, along with their elementary diagram, and no node may contain an element of the transitive closure of $p$ not given by $Z$. (This way, in the model induced by any infinite path through the tree, the interpretation of the symbol $``p"$ will be $p$ itself.) Furthermore, for a node $\tau$ to be a child of a node $\sigma$, all of the Skolem functions must have a value at $\tau$ on any input from $\sigma$. (This way, the model induced by any infinite path through the tree will be a model of $\varphi(p)$.) Conversely, any countable model of $\varphi(p)$ of height at most $\gamma$, along with a counting of that model, induces an infinite path through this tree. (To make choices among the successor of a node, you will in general need a choice function on this tree, which can be assumed to exist, since the tree can be constructed in $\rmL$ of some real parameters.)

Note that, for any $\gamma$, this tree is in $A$. Also, if it is well-founded, then it has a ranking function in any admissible set in which it is countable; and if it is not well-founded, then it has a branch definable over any such admissible set. It bears mention that there is such an admissible set in $A$, namely $\rmL_{\alpha}[Z]$, where $\alpha$ is the least $Z$-admissible beyond $\gamma$. So if $\varphi(p)$ is true in $V$, then there is a countable $\gamma$ such that the tree built on $\gamma$ is not well-founded, and hence has a branch definable over $\rmL_\alpha[Z]$, and therefore in $A$, from which a witness to $\varphi(p)$ can be built.
\end{proof}

Observe that by the upwards persistence of $\Sigma_1(p)$ formulas, if $\varphi$ is a $\Sigma_1(p)$ sentence true in $A$, then it is true in all supersets of $A$. Hence the absoluteness result Theorem~\ref{SL-lemma} also holds for all supersets of a transitive admissible set containing the countable ordinals in which $p$ is hereditarily countable.

Often the parameter in the statement of Theorem~\ref{SL-lemma} is taken to be a real $Z \subseteq \Nats$, in which case the statement can be simplified to absoluteness between $V$ and $\rmL_{\omega_1^{\rmL[Z]}}[Z]$. The formulation above is more general in that no coding $Z \subseteq \Nats$ of the transitive closure of $p$ need be assumed specified, there being no canonical choice of such a $Z$. The formulation above is superficially weaker in that, once $Z$ is chosen, $A$ must be at least $\rmL_{\omega_1}[Z]$, whereas in truth the potentially smaller set $\rmL_{\omega_1^{\rmL[Z]}}[Z]$ would suffice. However, this stronger version follows from the formulation given, since the theorem as stated gives absoluteness between $V$ and $\rmL[Z]$, and then the theorem could be interpreted in $\rmL[Z]$.


\begin{prop}
\label{Absolutness_of_asymp+}
	Let $L_0\subseteq L_1$ be countable languages, let $\M_0$ be an $L_0$-structure, and let $\M_1$ be an $L_1$-structure whose reduct to $L_0$ is $\M_0$. Let $A$ be a transitive admissible set containing the countable ordinals and in which $\M_0$ and $\M_1$ are hereditarily countable.
	Then the statement
	$\<e\>^{\M_0} \asymp^+ \M_1$ is absolute between $V$ and all supersets of $A$.
\end{prop}

\begin{proof}
Consider the definition of $\asymp^+$, Definition~\ref{fc-expansions-of-structures}. It is of the form ``for all $\M_0^*$ and $f$ and $\M_1^*$, where the $\M_j^*$'s are structures on $\Nats$, if $e$ feedback computes $\M_1^*$ from $\M_0^*$ then $\M_1^*$ is an expansion of $\M_0^*$ and the isomorphism $f$ extends to $\M_1^*$''. Whether $e$ feedback computes $\M_1^*$ from $\M_0^*$ is, by Lemma~\ref{remark-on-fb-natural-encoding}, $\Delta_1$ expressible as ``all admissible sets containing the $\M_j^*$'s satisfy a certain sentence" and ``there is an admissible set containing the $\M_j^*$'s satisfying a certain sentence". What follows is arithmetic in the parameters. Hence the entire definition is $\Pi_1$. By Theorem~\ref{SL-lemma}, applied to the negation of this relation, it is absolute between $V$ and $A$, and hence also with respect to all supersets of $A$.
\end{proof}

It follows from Proposition~\ref{Absolutness_of_asymp+} that the following definition is well-defined and doesn't depend on the specific forcing extension. Note that this can be seen as a feedback computability analogue of relations being uniformly relatively intrinsically (u.r.i.) computable (see \cite{AntonioBook}).

\begin{defi}
Let $L_0 \subseteq L_1$ be languages, $\M_0$ be a (not necessarily countable) $L_0$-structure, and $\M_1$ an $L_1$-structure which is an expansion of $\M_0$, i.e., $\M_1|_{L_0} = \M_0$. Then
	\defn{$e$ feedback computes the expansion $\M_1$ of oracle $\M_0$}, written
	$\<e\>^{\M_0} \asymp^+ \M_1$, if there is some forcing extension $\rmV[G]$ of the universe in which $\M_0$ is countable and $\rmV[G] \models \<e\>^{\M_0} \asymp^+ \M_1$.
\end{defi}

We now define what it means for a subset of a structure to be feedback computable.
\begin{defi}
Suppose $L_0$ is a language and $\M$ is a (not necessarily countable) $L_0$-structure. Suppose $U \subseteq \M$ is fixed by all automorphisms of $\M$ and $\M_U$ is the expansion of $\M$ which adds $U$ as a new unary predicate. Then $U$ is \defn{feedback computable} from $\M$ if there is an $e \in \Nats$ such that $\<e\>^{\M} \asymp^+ \M_U$.
\end{defi}

Note that by Proposition~\ref{Absolutness_of_asymp+} the specific forcing extension is irrelevant. The reason why we require $U$ to be closed under automorphisms of $\M$ is so that the set $U$ is uniquely defined by $\M_U$. Our main use of this notion is when $\M$ is of the form $(\gamma, \in_\gamma, A)$ for some ordinal $\gamma$ and finite subset $A$, where $\in_\gamma$ denotes the relation $\in$ restricted to $\gamma$.


\subsection{Feedback computing a structure}

Having defined what it means for an expansion of a structure to be feedback computable, we now define what it means for a structure to be feedback computable from another structure.

\begin{defi}
	\label{fc-between-structures}
Let $L_j$ be a language
	and $\M_j$ a countable $L_j$-structure for $j \in \{0, 1\}$. Then
	\defn{$e$ feedback computes $\M_1$ from $\M_0$}, written
	$\<e\>^{\M_0} \asymp \M_1$,
	if
	for every $L_0$-structure $\M_0^* \cong \M_0$ with underlying set $\Nats$ there is an $L_1$-structure $\M_1^* \cong \M_1$ with underlying set $\Nats$ such that
	$e$ feedback computes $\M_1^*$ from $\M_0^*$ using the natural encodings.
\end{defi}

As an example, let $G$ be a group and $H$ a normal subgroup. Consider the structure given by the group $G$ along with a distinguished relation for $H$. We can feedback compute the group $G/H$ by simply choosing a representative from each coset of $H$ along with the group multiplication table for these representatives induced by multiplying the corresponding cosets.

Note that for a structure $\M$ with underlying set $\Nats$ and natural encoding $i$, the notation $\<e\>^{i(\M)}$ denotes the feedback computable (partial) function from $\Nats$ to $\Nats$ that takes as an oracle the natural encoding of $\M$. In contrast, $\<e\>^\M$ will not be used on its own,
and $\<e\>^\M \asymp \N$ (or $\<e\>^\M \asymp^+ \N$) can be thought of as saying that no matter what copy of $\M$ is passed as an oracle to $\<e\>$, the output is always a copy of $\N$ (or in the case of $\asymp^+$ a copy of $\N$ which is also an expansion of $\M$).

It is worth noting that if $\<e\>^\M \asymp^+\N$ then we also have $\<e\>^\M \asymp \N$. However the converse need not hold as if $\<e\>^\M \asymp \N$ the output may be a structure whose restriction to $L_\M$ (the language of $\M$) is only isomorphic to $\M$ and not equal to it.

\begin{lem}
\label{lemma27}
Let $A$ be an admissible set in which models $\M$ and $\N$ (for the same language) are hereditarily countable. Then whether $\M$ and $\N$ are isomorphic in $V$ is $\Pi_1$ definable over $A$, and if they are isomorphic then an isomorphism is definable over $A$.
\end{lem}

\begin{proof}
Observe
that the set of isomorphisms is $\Sigma^1_1$ definable (in a real parameter coding $\M$ and $\N$ by any fixed standard way of coding a model by a real). Hence by the Kleene Basis Theorem, this set has a member computable in the hyperjump of the parameters (if non-empty) \cite[III.1.3]{HRT}.
\end{proof}

At its core, the proof is essentially building, in
$A$, the tree of finite partial isomorphisms between $\M$ and $\N$ (that is, isomorphisms between finite subsets of $\M$ and $\N$). They are isomorphic (in $V$) iff the tree is ill-founded iff there is no ranking function from the tree to the ordinals in $A$, a $\Pi_1$ statement over $A$; furthermore, an isomorphism can be built from any path through the tree, one of which is definable over $A$ (when the tree is ill-founded).

\begin{prop}
\label{Calculating-one-structure-from-another-absolute}
	Let $L_0$ and $L_1$ be languages, and let $\M_0$ be an $L_0$-structure and $\M_1$ an $L_1$-structure. Let $A$ be a transitive admissible set containing the countable ordinals and in which $\M_0$ and $\M_1$ are hereditarily countable. Then the statement $\<e\>^{\M_0} \asymp \M_1$ is absolute between $V$ and all supersets of $A$.
\end{prop}

\begin{proof}
Consider the definition of $\asymp$, Definition~\ref{fc-between-structures}. It is of the form ``for all $\M_0^*$ isomorphic to $\M_0$ there is an $\M_1^*$ such that $e$ feedback computes $\M_1^*$ from $\M_0^*$ and $\M_1^*$ and $\M_1$ are isomorphic''. As above, whether $e$ feedback computes $\M_1^*$ from $\M_0^*$ is, by Lemma~\ref{remark-on-fb-natural-encoding}, $\Delta_1$ expressible as ``all admissible sets containing the $\M_j^*$'s satisfy a certain sentence" and ``some admissible set containing the $\M_j^*$'s satisfies a certain sentence". Also, by Lemma~\ref{lemma27}, so is the assertion ``$\M_1^*$ and $\M_1$ are isomorphic", as a definable assertion over any admissible set in which the parameters are hereditarily countable. It bears observation that there is such an admissible set in $A$; namely, if $R \in A$ is a real coding (in some simple, standard way) both $\M_1$ and $\M_1^*$, then $\rmL_{\omega_1^R}[R] \in A$, as $A$ contains all countable ordinals and $\rmL_{\omega_1^R}[R] \subseteq A$. This provides a $\Pi_1$ definition of $\asymp$. By Theorem~\ref{SL-lemma}, applied to the negation of this relation, it is absolute between $V$ and $A$, and hence also with respect to all supersets of $A$.
\end{proof}

The absoluteness of the relation of feedback reducibility between structures
	(Definition~\ref{fc-between-structures})
allows us to make sense of feedback reducibility between structures even when those structures happen to be uncountable by considering the question of reducibility in a forcing extension where the structures are countable.

It follows from Proposition~\ref{Calculating-one-structure-from-another-absolute} that the following definition is well-defined and doesn't depend on the specific forcing extension.

\begin{defi}
Let $L_0, L_1$ be languages, $\M_0$ be a (not necessarily countable) $L_0$-structure, and $\M_1$ a (not necessarily countable) $L_1$-structure. Then
	\defn{$e$ feedback computes $\M_1$ from $\M_0$}, written
	$\<e\>^{\M_0} \asymp \M_1$, if there is some forcing extension $\rmV[G]$ of the universe in which $\M_0$ is countable and $\rmV[G] \models \<e\>^{\M_0} \asymp \M_1$.
\end{defi}

The relative computability of uncountable structures was studied using
generic extensions and Muchnik degrees
in \cite{MR3569106}. Our consideration of the feedback reducibility of uncountable structures can be seen as a feedback analogue of these notions, except using the analogue of Medvedev degrees instead of Muchnik degrees, because of the uniformity of our reductions.

Furthermore, both notions of feedback reducibility that take a structure as an oracle
(Definitions \ref{fc-expansions-of-structures} and \ref{fc-between-structures})
allow us to perform computation in a way that ignores the particular instantiations of the structures. This is important, as there are times when there is more computable information that can be obtained by the encoding of the structure than can be obtained intrinsically from the structure.

\subsection{Example: Functions from $\w$ to $\w$}

There is one example of computing one structure from another which is particularly important. Note that if $(W, {\lhd,}\,c)$ is a well-ordering of order type $\w$ and $c$ is a constant in $W$, then, uniformly in an encoding of an oracle $(W, \lhd, c)$, we can return the element of $\Nats$ which $c$ represents. Therefore if $\<e\>^{\M} \asymp (W, \lhd, c)$ there is little harm in identifying the output with the number $c$ represents. In particular we can define $\<e\>^{\M}(n) = m$ if $\<e\>^{\M \times (\w, \in, n)} \asymp (\w, \in, m)$. Hence we can think of $\<e\>^\M(n) = m$ as saying that whenever $e$ is handed a copy of $\M$ as an oracle, along with the natural number $n$, it outputs the natural number $m$.

As an example of this, suppose $G$ is a torsion group. There is a feedback computable function $\<e\>$ which takes $n$, computes the $n$th prime $p$ and returns the smallest $m > 0$ such that there is a subgroup of size $p^m$ if such an $m$ exists and returns $0$ otherwise. While the input depends on the specific group, it does not depend on the encoding of the group.

The following is then immediate from Proposition~\ref{Calculating-one-structure-from-another-absolute}.

\begin{cor}
Let $\M$ be an $L$-structure. The statement $\<e\>^{\M}(n) = m$
	is absolute among all transitive admissible sets
	containing the countable ordinals and in which $\M$ is hereditarily countable.
\end{cor}

\subsection{Example: $\alpha$-infinite time Turing machines}
For an admissible ordinal $\alpha$, the $(\alpha, \alpha)$-infinite time Turing machines (ITTMs) provide a different model which captures the $\alpha$-computable sets. We now show how to represent this model via oracle feedback computation, which highlights the value of feedback computation as a machine model for Borel maps.

The following definition is a straightforward generalization of the $(\infty, \w)$-Turing machines of \cite{HL} and the $(\alpha, \alpha)$-Turing machines in \cite{KS1}.

\begin{defi}
Let $\alpha, \beta$ be ordinals. A (run of a) \defn{$\alpha$-time, $\beta$-space Turing machine} (referred to as an $(\alpha, \beta)$-$\ITTM$) consists of the following data.
\begin{itemize}
\item A function $T\colon \alpha \times \beta \to \{0, 1\}$, called the \defn{tape}. For $\gamma \in \alpha$, the function $T(\gamma, \pars)\colon \beta \to \{0, 1\}$ is called the \defn{values} of the tape at time $\gamma$. The values at time $0$ are called the \defn{initial values}.

\item A function $H\colon \alpha \to \beta$, called the \defn{head location}.

\item A function $S\colon \alpha \to E$, called the \defn{state space}, where $E$ is a finite linearly ordered set containing a special starting state $s$ and halting state $h$. For $\gamma \in \alpha$, the value $S(\gamma)$ is called the \defn{state} of the machine at time $\gamma$.

\item A function $C\colon E \times \{0, 1\} \to E\times \{0, 1\} \times \{\LEFT, \RIGHT, \STAY\}$, called the \defn{lookup table}. It can be thought of as taking the state of the machine and the symbol written under the tape and returning the new state, the new symbol, and whether to move the head left or right, or to have it stay where it is.
\end{itemize}

This data is required to satisfy the following conditions.
\begin{itemize}
	\item $C(h, z)  = (h, z, \STAY)$ for all $z \in \{0,1\}$.

\item $H(0) = 0$ and $S(0) = s$.

\item If $\gamma \in \alpha$ is a limit ordinal then $H(\gamma)  = \liminf_{\zeta \in \gamma} H(\zeta)$ and $S(\gamma)  = \liminf_{\zeta \in \gamma} S(\zeta)$.

\item If $\gamma + 1 \in \alpha$ and $(e, z, M) = C\bigl(S(\gamma), T(\gamma, H(\gamma))\bigr)$, then the following hold.

\begin{itemize}
\item $S(\gamma + 1) = e$.

\item $T(\gamma+1, H(\gamma)) = z$.

\item $T(\gamma+1, \eta)  = T(\gamma, \eta)$ for $\eta \neq H(\gamma)$.

\item If $M = \STAY$ then $H(\gamma+ 1)  = H(\gamma)$.

\item If $M = \RIGHT$ then $H(\gamma+1)  = H(\gamma) + 1$.

\item If $M = \LEFT$ and $H(\gamma)  =  p + 1$ then $H(\gamma +1)  = p$.

\item If $M = \LEFT$ and $H(\gamma)$ is a limit ordinal then $H(\gamma +1)  = 0$.
\end{itemize}
\end{itemize}

	The \defn{input} of the machine is $T(0, \pars)$.
	The machine \defn{halts} if there is some $\gamma < \alpha$ such that $S(\gamma) = h$, and in this case, $T(\gamma, \pars)$ is the \defn{output} of the machine (which is well-defined by the first condition).

We will refer to an $(\alpha, \alpha)$-$\ITTM$ as simply an $\alpha$-$\ITTM$.
\end{defi}

We now show how to perform such computations using feedback.

\begin{lem}
\label{ITTM_feedback_computable_from_initial_conditions_and_code}
There is a feedback machine $\ittmcode$ such that if
	\[
		\M_{\alpha, \beta, X, C} \defas \<(\alpha, \in_\alpha), (\beta, \in_\beta), X, C\>
\]
		where $X\colon \beta\to \{0, 1\}$ and $C$ is a lookup table, then
\[
\<\ittmcode\>^{\M_{\alpha, \beta, X, C}} \cong \N_{\alpha, \beta, X, C}
\]
where $\N_{\alpha, \beta, X, C} = \<(\alpha, \in_\alpha), (\beta, \in_\beta), T_X, H_X, S_X, C\>$, with $(T_X, H_X, S_X, C)$ the (unique) $(\alpha, \beta)$-ITTM with $T_X(0, \pars) = X$ and code $C$.
\end{lem}
\begin{proof}
	Note that given $(\alpha, \beta)$, a code $C$, and initial values $X$,
	the definition of an $(\alpha, \beta)$-ITTM uniquely determines
	the functions $T_X$, $H_X$, and $S_X$
	by a transfinite recursion along $\alpha$ that is uniform in $\alpha$, $\beta$, and $X$.
Given a representation of an ordinal, we can feedback computably identify if an element of that representation corresponds to a limit ordinal (and if not find its successor).
	Hence from any isomorphic copy of $\M_{\alpha, \beta, X, C}$ we can feedback compute $\N_{\alpha, \beta, X, C}$.
\end{proof}

Henceforth all $(\alpha, \beta)$-ITTMs will have $\alpha = \beta$. As is standard in this situation, we will imagine that there is an input tape, an output tape, and an extra parameter tape (in which all but finitely many values are 0). This can be encoded into the ITTM in the standard way by interleaving these three tapes.

\begin{defi}
	A function $f \colon \Pfin(\alpha) \to \Pfin(\alpha)$ is \defn{$\alpha$-ITTM computable} if there is an $\alpha$-$\ITTM$ with a fixed finite extra parameter set such that when the input tape is the characteristic function of $A$ for some finite sequence $A$ of elements of $\alpha$, then the $\alpha$-$\ITTM$ halts with the characteristic function of $f(A)$ on the output tape.
The notion of $\alpha$-$\ITTM$ computability naturally extends to functions $f \colon \alpha^n \to \alpha^m$ for $n, m \in \Nats$.
\end{defi}

\begin{lem}
\label{ITTM_computable_ordinals_are_feedback_computable}
Let $\alpha$ and $\gamma$ be ordinals.
Suppose $(\alpha, \lhd)$ is well-ordered with order type $\gamma$ where $\lhd$ is
	$\alpha$-$\ITTM$ computable. Then for all finite $B \subseteq \gamma$ there is a feedback machine $e$ such that $\<e\>^{(\alpha, \in_\alpha, A)} \asymp (\gamma, \in_\gamma, B)$ for some finite $A \subseteq \alpha$.
\end{lem}
\begin{proof}
	For any finite lookup table $C$ and finite $A^*\in \Pfin(\alpha)$,
	from $(\alpha, \in_\alpha, A^*)$
	we can feedback compute $\M_{\alpha, \alpha, A^*, C}$ (via a feedback machine that intrinsically encodes $C$).
	Hence
	for any such $C$ and $A^*$, the structure
	$\N_{\alpha, \alpha, A^*, C}$ is feedback computable
	from $(\alpha, \in_\alpha, A^*)$
	by Lemma~\ref{ITTM_feedback_computable_from_initial_conditions_and_code}.
	By assumption, there is some $C$
	and some finite $A \in \Pfin(\alpha)$
	such that
from $\N_{\alpha, \alpha, A, C}$ we can feedback compute $(\gamma, \in_\gamma, B)$.
	Hence for some $A \in \Pfin(\alpha)$, the structure $(\gamma, \in_\gamma, B)$ can be feedback computed from $(\alpha, \in_\alpha, A)$.
\end{proof}

\subsection{Feedback computation from ordinals}

For the remainder of this section, we consider an extended example, feedback computability relative to a countable admissible ordinal.

\begin{prop}
\label{Feedback_in_alpha_implies_in_L(alpha+)}
Let $\gamma$ be an ordinal,
	let $A$ be a finite subset of $\gamma$, let $U\subseteq \gamma$, and suppose that $e$ is such that the function $\<e\>^{(\gamma, \in_\gamma, A)}$ feedback computes $(\gamma, \in_\gamma, U)$.
	Then $U \in L(\gamma^+)$.
\end{prop}
\begin{proof}
First note we can assume without loss of generality that $\gamma$ is countable, as if it isn't we can move to a forcing extension where $\gamma$ is countable.
Next note that the partial ordering $(\gamma^{<\w}, \preceq) \in L(\gamma^+)$ where $a \preceq b$ if $a$ is an initial segment of $b$. Let $G$ be a generic for $(\gamma^{<\w}, \preceq)$. Then $G$ is a surjection from $\w$ onto $\gamma$. By \cite[Theorem~1]{MR1159138}, the set $L(\gamma^+)[G]$ is admissible.
	From the surjection $G$ it is easy to find a real $G^* \in 2^\w$ encoding $(\gamma, \in_\gamma)$ such that $G^* \in L(\gamma^+)[G]$.

	Now let $M \in L(\gamma^+)[G]$ be a real encoding the structure $(\gamma, \in_\gamma, A)$ and let $M_U$ be the corresponding real encoding $(\gamma, \in_\gamma, A, U)$. Then $M_U$ is feedback computable from $M$ (by assumption). Therefore, by \cite[Proposition~16]{AFL}, the structure $M_U$ is in any admissible set containing $M$ and in particular is in $L(\gamma^+)[G]$. In particular this implies that $U \in L(\gamma^+)[G]$. But as $G$ was an arbitrary generic for $(\gamma^{<\w}, \preceq)$, we must have that $U \in L(\gamma^+)$, as desired.
\end{proof}

In general, even for a countable admissible $\alpha$, an ordinal $\gamma < \alpha^+$,
and a finite subset $A$ of $\gamma$,
there are not necessarily reals in $L(\alpha^+)$ encoding
$(\alpha, \in_\alpha)$ or $(\gamma, \in_\gamma, A)$.
However, Proposition~\ref{Feedback_in_alpha_implies_in_L(alpha+)} still
yields an upper bound on how complicated a function can be that is feedback computable from an arbitrary ordinal $\gamma$.

It is an interesting open question to pin down exactly how complicated the sets feedback computable from $\alpha$ can be. This is a question that we completely answer when $\alpha$ is a Gandy ordinal. It is worth pointing out that in Proposition~\ref{Feedback_in_alpha_implies_in_L(alpha+)} we cannot simply absorb the finite set $A$ into the code of the program, as the specific natural numbers representing the elements of $A$ depend on the particular representation of $(\gamma, \in_\gamma)$.

\begin{prop}
\label{Feedback_compute_L(gamma)_from_gamma}
	There is an $e$ such that $\<e\>^{(\gamma, \in_\gamma)} \asymp (L(\gamma), \in_{L(\gamma)})$.
	Further, $e$ is independent of $\gamma$.
\end{prop}
\begin{proof}
Without loss of generality we can assume $\gamma$ is countable.
	For a first-order formula $\varphi$ in the language of set theory,
	an ordinal $\alpha \in_\gamma \gamma$,
	a tuple $A$ of ordinals in $\alpha$ of length one less than the number of free variables in $\varphi$,
	define $\cS_{\varphi, A, \alpha} \defas \{ x \in L(\alpha) \st L(\alpha) \models \varphi(x, A)\}$.

Note that from
	$(\gamma, \in_\gamma)$ we can feedback compute the set of all such  $\cS_{\varphi, A, \alpha}$.
Further we can feedback compute the relation
	$\cS_{\varphi_0, A_0, \alpha_0} \in_{L(\gamma)} \cS_{\varphi_1, A_1, \alpha_1}$
	by induction on $\alpha_0$ and $\alpha_1$. Next, by induction on $\alpha$ we can compute an equivalence relation $\equiv$ where
	$\cS_{\varphi_0, A_0, \alpha_0} \equiv \cS_{\varphi_1, A_1, \alpha_1}$
	if and only if they contain the same elements.
	
	Let $i$ be a surjection from first-order formulas in the language of set theory to $\Nats$.
	For each $\equiv$-class,
	choose as a distinguished representative the $\cS_{\varphi, A, \alpha}$ where $(i(\varphi), A, \alpha)$ is lexicographically minimal.
	Finally, observe that
	$(L(\gamma), \in_{L(\gamma)})$ is isomorphic to the resulting collection of representatives under $\equiv$.
\end{proof}

We then have the following corollary.

\begin{cor}
\label{Gandy are feedback Gandy}
If there is an $\alpha$-computable well-ordering of $\alpha$ of height $\gamma$, then for any finite $B \subseteq \gamma$ there is an $e \in \Nats$ and a finite subset $A \subseteq \alpha$ such that $\<e\>^{(\alpha, \in_\alpha, A)} \asymp (\gamma, \in_\gamma, B)$.
\end{cor}
\begin{proof}
	We can feedback compute $L(\alpha+1)$ along with a well-ordering $\sqsubseteq_{L(\alpha+1)}$
	from $(\alpha, \in_\alpha)$, and so by Proposition~\ref{Feedback_compute_L(gamma)_from_gamma} we can feedback compute $L(\alpha+1)$ from $\alpha$. But every $\alpha$-computable well-ordering of $\alpha$ is in $L(\alpha+1)$,
	and so there must be some triple $(\varphi, A^*, \alpha)$ such that $\cS_{\varphi, A^*, \alpha}$ is an $\alpha$-computable well-ordering of $\alpha$ of order type $\gamma$.
	Therefore
	for some finite $A$ encoding $\varphi$, $A^*$, and $B$, we can compute $(\gamma, \in_\gamma, B)$.
\end{proof}

This suggests the following definition.

\begin{defi}
An admissible ordinal $\alpha$ is defined to be a \defn{feedback Gandy} ordinal if for all $\gamma < \alpha^+$ there is a well-ordering $(\alpha, \lhd)$ of order type $\gamma$ which is feedback computable from $(\alpha, \in_\alpha, A)$ for some finite $A \subseteq \alpha$.
\end{defi}

In particular, Corollary~\ref{Gandy are feedback Gandy} shows that all Gandy ordinals are feedback Gandy ordinals.
Whereas it has been established that there are admissible ordinals that are not Gandy ordinals \cite{MR552420}, it
is an open question whether or not every admissible ordinal is a feedback Gandy ordinal.

The following corollary is then immediate.

\begin{cor}
\label{Feedback_compute_from_alpha_ordinals_less_than_alpha+_for_Gandy_ordinals}
If $\alpha$ is a feedback Gandy ordinal and $\gamma < \alpha^+$, then for all finite $B \subseteq \gamma$ there is a finite $A \subseteq \alpha$ such that $(\gamma, \in_\gamma, B)$ is feedback computable from $(\alpha, \in_\alpha, A)$.
\end{cor}

\begin{prop}
\label{Elements in L(gamma) are feedback computable from alpha when gamma is}
Let $\alpha < \gamma$ be ordinals, let $A$ be a finite subset of $\alpha$, and let $e \in \Nats$. Suppose that
	$U \in L(\gamma)$ is such that $U \subseteq \alpha$,
	and suppose that
	$\<e\>^{(\alpha, \in_\alpha, A)}\asymp (\gamma, \in_\gamma)$. Then there is some finite $A^* \subseteq \alpha$ and $e^* \in \Nats$ such that $\<e^*\>^{(\alpha, \in_\alpha, A^*)} \asymp^+ (\alpha, \in_\alpha, U)$.
\end{prop}
\begin{proof}
This follows immediately from Proposition~\ref{Feedback_compute_L(gamma)_from_gamma} by letting $A^* \defas A \cup \{a\}$, where $a$ is a code for $U$ in a definable bijection from $\gamma$ to $L(\gamma)$.
\end{proof}

Combining Propositions~\ref{Feedback_in_alpha_implies_in_L(alpha+)} and
\ref{Elements in L(gamma) are feedback computable from alpha when gamma is}
and Corollary~\ref{Feedback_compute_from_alpha_ordinals_less_than_alpha+_for_Gandy_ordinals},
we obtain the following.

\begin{thm}
If $\alpha$ is a feedback Gandy ordinal then the following are equivalent for $U \subseteq \alpha$.

\begin{itemize}
\item $U \in L(\alpha^+)$.

\item $U$ is feedback computable from $(\alpha, \in_\alpha, A)$ for some $A \in \Pfin(\alpha)$.
\end{itemize}
\end{thm}

\section{Open questions}

We end with several open questions. For each of these questions, let $\alpha$ be an ordinal and $A$ a finite subset of $\alpha$.

\begin{itemize}
\item For what $\beta$ is there some set $U \in L(\beta+1)\setminus L(\beta)$ that is feedback computable from $(\alpha, \in_\alpha, A)$?

\item If there is some $U \in L(\beta+1)\setminus L(\beta)$ that is feedback computable from $(\alpha, \in_\alpha, A)$, must $L(\beta+1)$ be feedback computable from $(\alpha, \in_\alpha, A^*)$ for some finite $A^* \subseteq \alpha$?

\item Which ordinals are feedback Gandy? In particular, are there feedback Gandy ordinals that are not Gandy ordinals? Indeed, are there any admissible ordinals that are not feedback Gandy ordinals?

\end{itemize}

\section*{Acknowledgement}
The authors would like to thank Julia Knight, Russell Miller, Noah Schweber, and Richard Shore for helpful conversations, and the anonymous referees for their useful comments.


\providecommand{\bysame}{\leavevmode\hbox to3em{\hrulefill}\thinspace}
\providecommand{\MR}{\relax\ifhmode\unskip\space\fi MR }
\providecommand{\MRhref}[2]{%
  \href{http://www.ams.org/mathscinet-getitem?mr=#1}{#2}
}
\providecommand{\href}[2]{#2}


\begin{thebibliography}{KMS16}

\bibitem[AFL15]{AFL}
N.~L. Ackerman, C.~E. Freer, and R.~S. Lubarsky, \emph{Feedback {T}uring
  computability, and {T}uring computability as feedback}, 30th Annual
  {ACM/IEEE} Symposium on Logic in Computer Science, {LICS} 2015, Kyoto, Japan,
  {IEEE}, 2015, pp.~523--534.

\bibitem[AS76]{MR0432431}
F.~G. Abramson and G.~E. Sacks, \emph{Uncountable {G}andy ordinals}, J. London
  Math. Soc. (2) \textbf{14} (1976), no.~3, 387--392.

\bibitem[Ers90]{MR1159138}
Y.~L. Ershov, \emph{Forcing in admissible sets}, Algebra and Logic \textbf{29}
  (1990), no.~6, 424--430 (1991).

\bibitem[Gos79]{MR552420}
R.~Gostanian, \emph{The next admissible ordinal}, Ann. Math. Logic \textbf{17}
  (1979), no.~1-2, 171--203.

\bibitem[Hin12]{MR2931672}
P.~G. Hinman, \emph{A survey of {M}u\v cnik and {M}edvedev degrees}, Bull.
  Symbolic Logic \textbf{18} (2012), no.~2, 161--229.

\bibitem[HL00]{HL}
J.~D. Hamkins and A.~Lewis, \emph{Infinite time {T}uring machines}, J. Symbolic
  Logic \textbf{65} (2000), no.~2, 567--604.

\bibitem[Hoo82]{Hoo}
M.~R.~R. Hoole, \emph{Recursion on sets}, Ph.D. thesis, Univ. of Oxford, 1982.

\bibitem[Jec97]{MR1492987}
T.~Jech, \emph{Set theory}, 2nd ed., Perspectives in Mathematical Logic,
  Springer-Verlag, Berlin, 1997.

\bibitem[Kec95]{Kechris}
A.~S. Kechris, \emph{Classical descriptive set theory}, Graduate Texts in
  Mathematics, vol. 156, Springer-Verlag, New York, 1995.

\bibitem[KMS16]{MR3569106}
J.~Knight, A.~Montalb\'an, and N.~Schweber, \emph{Computable structures in
  generic extensions}, J. Symb. Log. \textbf{81} (2016), no.~3, 814--832.

\bibitem[KS09]{KS1}
P.~Koepke and B.~Seyfferth, \emph{Ordinal machines and admissible recursion
  theory}, Ann. Pure Appl. Logic \textbf{160} (2009), no.~3, 310--318.

\bibitem[Lub88]{Lub}
R.~S. Lubarsky, \emph{Another extension of {V}an de {W}iele's theorem}, Ann.
  Pure Appl. Logic \textbf{38} (1988), no.~3, 301--306.

\bibitem[Mon18]{AntonioBook}
A.~Montalb\'{a}n, \emph{Computable structure theory}, Draft of Part I,
  \texttt{https://math.berkeley.edu/} \texttt{\textasciitilde antonio/CSTpart1.pdf}, 2018.

\bibitem[Mos09]{M}
Y.~N. Moschovakis, \emph{Descriptive set theory}, second ed., Mathematical
  Surveys and Monographs, vol. 155, American Mathematical Society, Providence,
  RI, 2009.

\bibitem[Sac90]{HRT}
G.~E. Sacks, \emph{Higher recursion theory}, Perspectives in Mathematical
  Logic, Springer-Verlag, Berlin, 1990.

\bibitem[Sla85]{Sl1}
T.~A. Slaman, \emph{Reflection and forcing in {$E$}-recursion theory}, Ann.
  Pure Appl. Logic \textbf{29} (1985), no.~1, 79--106.

\bibitem[Sla86]{Sl2}
\bysame, \emph{{$\Sigma_1$}-definitions with parameters}, J. Symbolic Logic
  \textbf{51} (1986), no.~2, 453--461.

\bibitem[vdW82]{vdw}
J.~van~de Wiele, \emph{Recursive dilators and generalized recursions}, Proc.
  Herbrand Symposium, Logic Colloquium '81 (J.~Stern, ed.), North-Holland,
  Amsterdam, 1982, pp.~325--332.

\end{thebibliography}
\end{document}